\documentclass[12pt]{article}
\usepackage{amssymb,amsmath,amsthm,amsfonts}
\usepackage{epic,eepic}
\newtheorem{thm}{Theorem}[section]
\newtheorem{lemma}[thm]{Lemma}

\newcommand{\prlabel}[1]{\label{#1}}
\newcommand{\prbibitem}[2]{\bibitem[#1]{#2}}
\newcommand{\R}{{\mathbb{R}}}

\newcommand{\T}{{\mathbb{T}}}
\newcommand{\Z}{{\mathbb{Z}}}
\newcommand{\N}{{\mathbb{N}}}
\newcommand{\F}{{\mathbb{F}}}
\newcommand{\la}{{\lambda}}

\newcommand{\Ar}{{\text{Area}}}
\newcommand{\Le}{{\text{Length}}}
\newcommand{\QED}{{\square}}
\begin{document}

\title{A linear isoperimetric inequality for the punctured
Euclidean plane
}

\author{Leonid Polterovich\thanks{Supported by 
THE ISRAEL SCIENCE FOUNDATION founded by the Israel Academy
of Sciences and Humanities.}
\\
Tel Aviv
University, Israel
\\ polterov@post.tau.ac.il
\;\and Jean-Claude Sikorav\\
E.N.S. de Lyon, France
\\Jean-Claude.SIKORAV@umpa.ens-lyon.fr    }

\date{June 28, 2001}

\maketitle

\centerline{preliminary version}

\begin{abstract}

It follows from a general theorem of Bonk and Eremenko 
that contractible 
closed curves on the Euclidean $\R^2 \setminus \Z^2$
satisfy a linear isoperimetric inequality. 
In the present note we give an alternative proof of this
fact. Our approach is based on a non-standard
combinatorial isoperimetric inequality
for the group $<a,b,c| aba^{-1}b^{-1}c>$.
We show that the combinatorial area of every element in
the normal closure of $aba^{-1}b^{-1}c$ does not exceed 
the half of its length in generators $a^{\pm 1},b^{\pm 1}$
(here $c^{\pm 1}$ is ignored). The proof uses
a refinement of the small cancellation theory. In addition, 
we present an application of the isoperimetric inequality
on $\R^2 \setminus \Z^2$  
to Hamiltonian dynamics. Combining it
with methods of symplectic topology
we show that every non-identical
Hamiltonian diffeomorphism of the 2-torus has at least linear
asymptotic growth of the differential.

\end{abstract}

\newpage

\section{Introduction and results}

For a smooth contractible curve $\alpha :S^1 \to \R^2 \setminus \Z^2$
define its area by
$$\Ar (\alpha) = \inf \int_{D^2} |f^*\omega|,$$
where the infimum is taken over all maps 
$f: D^2 \to \R^2 \setminus \Z^2$
with $f|_{S^1} = \alpha$, and 
$\omega$ stands for the Euclidean area form.
A recent result due to Bonk and Eremenko \cite{BE}
which is proved for a much more general class of surfaces
yields that $\Ar (\alpha) \leq \mu \Le (\alpha)$
for some universal constant $\mu$.
The proof in \cite {BE} uses the theory of Gromov hyperbolic
spaces.
In this note we give an alternative proof of this result
and present an application to Hamiltonian dynamics.
\medskip

\begin{thm}[see Bonk-Eremenko\cite{BE}]
\prlabel{ineq1}
$$\Ar (\alpha) \leq (1 + \sqrt{2}) \Le (\alpha)$$
for every smooth contractible 
curve $\alpha :S^1 \to \R^2 \setminus \Z^2$.
\end{thm}

\medskip
\noindent

Theorem \ref{ineq1} has the following counterpart in
combinatorial group theory. Consider the free group
$\F_3$ with 3 generators $a,b,c$ and the word $r = aba^{-1}b^{-1}c$.
Let $N \subset F$ be the normal closure of $r$. For an element
$w \in N$ define its {\it combinatorial area} $A(w) = \inf d$, where
the infimum is taken over all decompositions of the form
$$ w = \prod_{i=1}^d g_i r^{\pm 1} g_i ^{-1}, \;\; g_i \in \F_3.$$
Define the {\it $(a,b)$ -length} $l(w)$ as the number 
of occurrences of symbols
$a^{\pm 1}, b^{\pm 1}$ in a cyclic reduction of $w$.

\medskip

\begin{thm}
\prlabel{ineq2}
$A(w) \leq \frac{1}{2} l(w)$ for every $w \in N$.
\end{thm}

\medskip
\noindent
The link between these two inequalities is discussed in \S 2
where we deduce Theorem \ref{ineq1} from Theorem \ref{ineq2}.
Theorem \ref{ineq2} is proved in \S 3.

Some remarks are in order.
Inequality \ref{ineq2} is sharp (see 3.3 below)
while the isoperimetric constant provided by
Theorem \ref{ineq1} most probably can be improved. It would 
be interesting to compare the isoperimetric
constant from \ref{ineq1} with the one obtained in \cite{BE}.

For an embedded curve $\alpha$ one has $\Ar (\alpha) \leq \Le (\alpha)$.
This follows from the classical theorem of Jarnik in number theory
(see \cite{Hua}, p. 123).

Contractibility of $\alpha$ cannot be replaced by a homological
assumption. Indeed,
let $\la$ be a primitive of $\omega$ on $\R^2$, that is $\omega = d\la$.
It is easy to construct a curve $\alpha: S^1 \to \R^2 \setminus \Z^2$
of arbitrarily large length which 
is homologous to zero (but non-contractible!)
so that
$$\int_{\alpha} \la \geq {\Big (} \Le (\alpha) {\Big )}^ {1+ \delta}
$$
for some positive constant $\delta$.

A linear isoperimetric inequality is obvious for
the complement of the ``thick" lattice.
Let $B(\epsilon)$ be the
Euclidean $\epsilon$-neighbourhood
of $\Z^2$ in $\R^2$.  One can easily show  that
$\Ar (\alpha) \leq C(\epsilon) \Le (\alpha)$ for every 
contractible curve
$\alpha : S^1 \to \R^2 \setminus B(\epsilon)$. Indeed, the Euclidean
metric on $\R^2 \setminus B(\epsilon)$
is equivalent to the hyperbolic one
which as it is known satisfies the
linear isoperimetric inequality. The problem is that the
Lipschitz constant, and hence the isoperimetric constant
obtained by this argument, goes to infinity as $\epsilon \to 0$.
Bruce Kleiner pointed out that one can go round this difficulty
using the theory of Gromov hyperbolic spaces. He outlined an
argument which shows that
a linear isoperimetric inequality should be valid
for much more
general spaces of non-positive curvature . 
Bonk and Eremenko \cite{BE}
established a linear isoperimetric inequality for a
fairly general class of open (not necessarily complete)
surfaces of non-positive curvature.

The inequality $A(w) \leq  |w|$
where $|w|$ stands for the usual word length of $w \in N$
in generators $a,b,c$ is an immediate consequence of the
small cancellation theory \cite{LS}. It turns out that in our situation
some arguments of this theory can be refined. This enables us to
replace the word length
by the $(a,b)$-length  and to improve the isoperimetric
constant.

Finally let us discuss an application of isoperimetric inequality
\ref{ineq1} to Hamiltonian dynamics.
Consider an area-preserving diffeomorphism $\psi: \R^2 \to \R^2$
of the form
$$\psi(x_1,x_2) = (x_1+ \psi_1(x_1,x_2),x_2 + \psi_2(x_1,x_2)),$$
where the functions $\psi_i,\;i=1,2$ are $1$-periodic in both variables
and satisfy
\begin{equation}
\prlabel{ham}
\int_0^1\int_0^1 \psi_i(x_1,x_2)dx_1dx_2
= 0.
 \end{equation}

Diffeomorphisms of the torus $\T^2 = \R^2/\Z^2$ corresponding
to such $\psi$'s are called {\it Hamiltonian}. They play a fundamental
role both in symplectic topology and classical mechanics
(see \cite{MS},\cite{HZ},\cite{P1}).
Put
$$||d\psi|| = \max |d_x\psi (\xi)|,$$
where the maximum is taken over all unit 
tangent vectors $\xi \in T_x \R^2$.
The growth type of the sequence $||d\psi^n||$ is a basic 
dynamical invariant
of $\psi$ (cf. \cite{DG}). Up to a multiplicative constant 
it is an invariant of
the conjugacy class of $\psi$ in the group 
of all $\Z^2$-periodic diffeomorphisms.

\medskip
\noindent
\begin{thm}
\prlabel{growth}
Suppose that $\psi \neq \text{id}$. Then
$$||d\psi^n || \geq \kappa n$$
for some $\kappa = \kappa(\psi) > 0$.
\end{thm}

\medskip
\noindent
This result is sharp. For instance $||d\psi^n||$ grows linearly
for $\psi(x_1,x_2) = (x_1,x_2 + \sin x_1)$. On the other hand
one has $||d\phi^n|| \equiv 1$ when $\phi$ is a parallel
translation $x \to x+ \text{const}$ of the plane.
Here $\phi$ is still area-preserving and $\Z^2$-periodic,
but non-Hamiltonian -- condition (\ref{ham}) is violated.
We refer to \cite{P2} for
an extensive discussion on the growth of differential
of symplectic maps including
more sophisticated non-Hamiltonian  examples and
counterexamples on $\T^2$,
sharp results for surfaces of higher genus and
generalizations to Hamiltonian diffeomorphisms
of arbitrary aspherical symplectic manifolds.
Theorem \ref{growth} is proved in \S 4 below.


\section{Reduction to an algebraic problem}

Let $G_0$ be the unit grid on the plane with integer vertices.
Fix $\epsilon \in (0;0.1)$
and write $B(\epsilon)$ for the $\epsilon$-neighbourhood of
$\Z^2$ in $\R^2$. Consider the 1-complex
$$G_{\epsilon} = (G_0 \setminus B(\epsilon)) \cup \partial B(\epsilon)$$
(see figure 1). A 1-cycle on 
$G_{\epsilon}$ is called {\it geometrically irreducible}
if it has no pair of consecutive edges which coincide but have opposite
orientations. We use the convention that a trivial 1-cycle
(a vertex) is geometrically irreducible.
Clearly, every closed curve on $G_{\epsilon}$ is free homotopic
to a geometrically irreducible cycle, and
every
 geometrically irreducible cycle minimizes the Euclidean
length among closed curves  on $G_{\epsilon}$
in its free homotopy class. This follows from the fact 
that the universal
cover of $G_{\epsilon}$ is a tree.

\begin{figure}[h]
\label{fig1}
$$\setlength{\unitlength}{0.00083333in}
\begingroup\makeatletter\ifx\SetFigFont\undefined%
\gdef\SetFigFont#1#2#3#4#5{%
  \reset@font\fontsize{#1}{#2pt}%
  \fontfamily{#3}\fontseries{#4}\fontshape{#5}%
  \selectfont}%
\fi\endgroup%
{\renewcommand{\dashlinestretch}{30}
\begin{picture}(3324,3339)(0,-10)
\put(612,2712){\ellipse{600}{600}}
\put(2721,2710){\ellipse{600}{600}}
\put(2721,610){\ellipse{600}{600}}
\put(612,612){\ellipse{600}{600}}
\put(2712,312){\blacken\ellipse{90}{90}}
\put(2712,312){\ellipse{90}{90}}
\put(2412,612){\blacken\ellipse{90}{90}}
\put(2412,612){\ellipse{90}{90}}
\put(912,612){\blacken\ellipse{90}{90}}
\put(912,612){\ellipse{90}{90}}
\put(312,612){\blacken\ellipse{90}{90}}
\put(312,612){\ellipse{90}{90}}
\put(612,912){\blacken\ellipse{90}{90}}
\put(612,912){\ellipse{90}{90}}
\put(612,312){\blacken\ellipse{90}{90}}
\put(612,312){\ellipse{90}{90}}
\put(612,2412){\blacken\ellipse{90}{90}}
\put(612,2412){\ellipse{90}{90}}
\path(612,3312)(612,3012)
\path(12,2712)(312,2712)
\path(12,612)(312,612)
\path(612,12)(612,312)
\path(2712,12)(2712,312)
\path(3312,612)(3012,612)
\path(3312,2712)(3012,2712)
\path(2712,3312)(2712,3012)(2712,3012)
\path(912,2712)(2412,2712)
\path(2712,2412)(2712,912)
\path(2412,612)(912,612)
\path(612,2412)(612,912)
\put(312,2187){\makebox(0,0)[lb]{$P_2$}}
\put(312,987){\makebox(0,0)[lb]{$M$}}
\put(87,387){\makebox(0,0)[lb]{$E_1$}}
\put(762,87){\makebox(0,0)[lb]{$P_1$}}
\put(987,687){\makebox(0,0)[lb]{$L$}}
\put(2187,387){\makebox(0,0)[lb]{$E_2$}}
\put(2862,87){\makebox(0,0)[lb]{$P_3$}}
\end{picture}
}
$$
\caption{}
\end{figure}
For a homotopy $\theta: S^1 \times [0;1] \to \R^2$ define its {\it area}
as the integral of $|\theta^*\omega|$ over $S^1 \times [0;1]$.

\medskip

\begin{lemma}
\prlabel{hom}
Every smooth curve $\alpha : S^1 \to \R^2 \setminus B(\epsilon)$
can be homotoped in $\R^2 \setminus \Z^2$  
to a geometrically irreducible cycle $\gamma : S^1 \to G_{\epsilon}$
with $\Le (\gamma) \leq 2\Le (\alpha)$
by a homotopy of area $\leq \sqrt{2} \Ar (\alpha)$.
\end{lemma}

\medskip
\noindent
The proof is given at the end of this section.
Let us now look more attentively at curves on $G_{\epsilon}$.
First consider the complex $G_{\epsilon}/\Z^2$.
Its fundamental group equals $\F_3$.
Let $a,b$
and $c$ be generators of this fundamental group represented
by the following paths on $G_{\epsilon}$ (see figure 1): $a = P_1E_1MP_2,
b= P_1LE_2P_3$ and $c = P_1LME_1P_1$. In Euclidean coordinates
we have
$P_1 = (0,-\epsilon)$, $P_2= (0,1-\epsilon)$, $P_3 = (1,-\epsilon)$,
$E_1 = (-\epsilon,0)$, $E_2 = (1-\epsilon, 0)$, $L = (\epsilon, 0)$
and $M = (0,\epsilon)$. Points denoted by the same letter with different
indices have the same projections to $G_{\epsilon}/\Z^2$. The fundamental
group described above corresponds to the base point $P_1 \; \mod \; \Z^2$.

Consider the homomorphism of $\F_3 = <a,b,c>$ to $\Z^2$
which sends $a$ and $b$ to $(0,1)$ and $(1,0)$ respectively
and vanishes on $c$.
Every word $w \in \F_3$ which lies in its kernel
represents a closed curve on $G_{\epsilon}$ denoted
$\gamma '(w)$.
By Seifert-van Kampen theorem this curve is contractible
on $\R^2 \setminus B(\epsilon)$ if and only if
$w$ belongs to the normal closure $N$ of the element
$r = aba^{-1}b^{-1}c$. Indeed, $\R^2 \setminus B(\epsilon)$
is a cover of  the surface which is obtained from
$G_{\epsilon} /\Z^2$ by gluing a disc along the cycle $r$.
Write $\gamma (w)$ for the geometrically irreducible cycle
homotopic to $\gamma ' (w)$.

\medskip
\noindent
\begin{lemma}
\prlabel{over}
Let $w \in N$ be a cyclically reduced word.
 Then
$$l(w) \leq (1-2\epsilon)^{-1} \Le (\gamma(w)).$$
\end{lemma}

\medskip
\noindent
{\bf Proof:} Note that consecutive edges of $\gamma '(w)$ which coincide
but have opposite orientation can occur only on circular parts of curves
corresponding to letters $a^{\pm 1}, b^{\pm 1}$ in $w$. The length
of the segment of the grid $G_0$ corresponding to each such letter
is $1-2\epsilon$. Hence $\Le (\gamma (w)) \geq (1-2\epsilon)l(w)$.
$\QED$

\medskip
\noindent
{\bf Proof of Theorem \ref{ineq1}:}
Let $\alpha$ be a contractible curve on $\R^2  \setminus \Z^2$.
Fix $\epsilon > 0$  so small that $\alpha$ lies in $\R^2 \setminus B(\epsilon)$.
Lemma \ref{hom}
provides a homotopy of $\alpha$ to a geometrically irreducible cycle
$\gamma$ on $G(\epsilon)$ with
$\Le (\gamma) \leq 2 \Le (\alpha)$.
Let $w \in N$ be a cyclically reduced word so that
$\gamma = \gamma (w)$. Clearly, $\Ar (\gamma) \leq (1-\pi \epsilon^2) A(w)$.
By Theorem \ref{ineq2} $A(w) \leq l(w)/2$. Combining this with Lemma \ref{over}
we get that
$$\Ar(\gamma) \leq (1-\pi \epsilon^2) (1-2\epsilon)^{-1} \Le (\alpha).$$
Recall now that
the area of the homotopy between $\alpha$ and $\gamma$ does not exceed
$\sqrt{2}\Le (\alpha)$. Hence, taking $\epsilon \to 0$ we get
$$\Ar (\alpha) \leq (1 + \sqrt{2} )\Le (\alpha)$$ as required.
$\QED$

\medskip
\noindent
{\bf Proof of Lemma \ref{hom}:}
 Given a curve
$\alpha: S^1 \to \R^2 \setminus B(\epsilon)$,
perturb it so that it will be in general position with respect
to $G_{\epsilon}$. This means that $G_{\epsilon}$ does not contain
self-intersection points of $\alpha$ and $\alpha$ intersects
$G_{\epsilon}$ transversally. 
Assume first that $\alpha \cap G_{\epsilon} \neq \emptyset$.
Then there exists a finite partition
$[0;1] = \bigcup_{i=1}^{n} I_i$ such that $\alpha |_{I_{i}}$
lies entirely in a fundamental domain of $\R^2  \setminus B(\epsilon)$.
Fix a segment $I = [t_i;t_{i+1}]$ of the 
partition, and put $X = \alpha (t_i),
Y = \alpha (t_{i+1})$ and $f = \alpha|_I$. Denote by $Q$ the fundamental
domain containing $f$.
There are 3 possibilities.

\medskip
\noindent
{\bf Case 1.} Points $X$ and $Y$ lie on different and non-opposite
linear edges of $Q$ (see figure 2a).

\medskip
\noindent
{\bf Case 2.} Points $X$ and $Y$ lie on opposite
linear edges of $Q$ (see figure 2b).

\medskip
\noindent
{\bf Case 3.} Points $X$ and $Y$ lie on the
same linear edge of $Q$ (see figure 2c).

\medskip
\noindent
Let us outline the proof of the Lemma.
In each case, we are going to homotope (with fixed end points)
the curve $f$ to a curve
$g$ which is a part of $\partial Q$ between $X$ and $Y$.
We will see that
$\Le (g) \leq 2 \Le (f)$ and the area of the homotopy does not
exceed $\sqrt{2} \Le (f)$. Clearly, this will enable us to homotope
$\alpha$ to a curve $\gamma $ on $G(\epsilon)$ which satisfies the
inequalities of the Lemma. Shortening $\gamma $ if necessary, we
achieve that $\gamma$ is geometrically irreducible.
We treat each of these three cases separately.

\begin{figure}[h]\label{fig2}
$$
\setlength{\unitlength}{0.00045833in}
\begingroup\makeatletter\ifx\SetFigFont\undefined%
\gdef\SetFigFont#1#2#3#4#5{%
  \reset@font\fontsize{#1}{#2pt}%
  \fontfamily{#3}\fontseries{#4}\fontshape{#5}%
  \selectfont}%
\fi\endgroup%
{\renewcommand{\dashlinestretch}{30}
\begin{picture}(10713,3550)(0,-10)
\put(187.500,2816.500){\arc{855.132}{0.2663}{1.3045}}
\put(2212.500,2816.500){\arc{855.132}{1.8370}{2.8753}}
\put(2212.500,791.500){\arc{855.132}{3.4078}{4.4461}}
\put(187.500,791.500){\arc{855.132}{4.9786}{6.0169}}
\path(300,2404)(300,1204)
\path(600,904)(1800,904)
\path(2100,1204)(2100,2404)
\path(1800,2704)(600,2704)
\put(4087.500,2816.500){\arc{855.132}{0.2663}{1.3045}}
\put(6112.500,2816.500){\arc{855.132}{1.8370}{2.8753}}
\put(6112.500,791.500){\arc{855.132}{3.4078}{4.4461}}
\put(4087.500,791.500){\arc{855.132}{4.9786}{6.0169}}
\path(4200,2404)(4200,1204)
\path(4500,904)(5700,904)
\path(6000,1204)(6000,2404)
\path(5700,2704)(4500,2704)
\put(7987.500,2816.500){\arc{855.132}{0.2663}{1.3045}}
\put(10012.500,2816.500){\arc{855.132}{1.8370}{2.8753}}
\put(10012.500,791.500){\arc{855.132}{3.4078}{4.4461}}
\put(7987.500,791.500){\arc{855.132}{4.9786}{6.0169}}
\path(8100,2404)(8100,1204)
\path(8400,904)(9600,904)
\path(9900,1204)(9900,2404)
\path(9600,2704)(8400,2704)
\Thicklines
\put(187.500,791.500){\arc{855.132}{4.9786}{6.0169}}
\put(4087.500,2816.500){\arc{855.132}{0.2663}{1.3045}}
\put(4087.500,791.500){\arc{855.132}{4.9786}{6.0169}}
\thinlines
\put(2325.000,1054.000){\arc{1209.339}{5.2315}{6.4075}}
\put(296,908){\blacken\ellipse{76}{76}}
\put(296,908){\ellipse{76}{76}}
\put(4196,908){\blacken\ellipse{76}{76}}
\put(4196,908){\ellipse{76}{76}}
\put(1355,149){\ellipse{282}{282}}
\put(5255,168){\ellipse{282}{282}}
\put(9155,168){\ellipse{282}{282}}
\path(300,904)(3000,904)
\path(2880.000,874.000)(3000.000,904.000)(2880.000,934.000)
\Thicklines
\path(300,2104)(300,1204)
\drawline(600,904)(600,904)
\path(600,904)(1350,904)
\path(4950,2704)(4500,2704)
\path(4200,2404)(4200,1204)
\path(4500,904)(4950,904)
\path(8100,2179)(8100,1504)
\thinlines
\path(4226.000,3188.000)(4196.000,3308.000)(4166.000,3188.000)
\path(4196,3308)(4196,908)
\path(4200,904)(6600,904)
\path(6480.000,874.000)(6600.000,904.000)(6480.000,934.000)
\path(8130.000,3184.000)(8100.000,3304.000)(8070.000,3184.000)
\path(8100,3304)(8100,904)
\path(8100,904)(10500,904)
\path(10380.000,874.000)(10500.000,904.000)(10380.000,934.000)
\path(2659,1489)(2603,1620)(2734,1564)
\path(300,2104)(303,2103)(310,2101)
	(322,2097)(338,2091)(358,2084)
	(380,2077)(402,2069)(424,2062)
	(445,2054)(464,2048)(482,2041)
	(500,2035)(516,2029)(533,2023)
	(550,2016)(566,2011)(582,2005)
	(598,1998)(616,1992)(634,1985)
	(653,1978)(672,1971)(691,1964)
	(711,1957)(730,1951)(749,1945)
	(768,1939)(786,1934)(803,1929)
	(819,1925)(834,1922)(849,1919)
	(863,1916)(879,1915)(895,1914)
	(910,1913)(925,1914)(939,1916)
	(951,1918)(963,1922)(973,1926)
	(981,1930)(988,1934)(994,1939)
	(997,1944)(999,1949)(1000,1954)
	(999,1960)(996,1965)(992,1971)
	(986,1977)(979,1982)(972,1986)
	(963,1990)(954,1993)(946,1994)
	(938,1995)(931,1994)(925,1992)
	(925,1991)(920,1988)(915,1983)
	(912,1976)(909,1967)(907,1957)
	(907,1945)(907,1932)(909,1918)
	(912,1903)(915,1887)(920,1871)
	(925,1854)(930,1841)(935,1827)
	(940,1812)(947,1797)(953,1780)
	(960,1763)(968,1745)(975,1727)
	(982,1709)(990,1691)(997,1673)
	(1003,1655)(1010,1639)(1015,1622)
	(1020,1607)(1025,1591)(1030,1572)
	(1035,1552)(1038,1532)(1041,1513)
	(1042,1494)(1043,1476)(1042,1459)
	(1041,1444)(1038,1432)(1034,1422)
	(1030,1414)(1024,1408)(1020,1405)
	(1016,1403)(1011,1402)(1006,1401)
	(1000,1402)(995,1403)(989,1405)
	(984,1407)(978,1411)(973,1414)
	(969,1419)(965,1424)(962,1429)
	(960,1434)(959,1439)(958,1445)
	(958,1450)(960,1456)(963,1464)
	(968,1471)(975,1479)(983,1487)
	(993,1495)(1004,1502)(1016,1509)
	(1028,1516)(1041,1521)(1054,1526)
	(1066,1529)(1078,1532)(1089,1533)
	(1099,1533)(1109,1531)(1118,1529)
	(1127,1525)(1135,1519)(1143,1513)
	(1150,1505)(1157,1496)(1163,1486)
	(1169,1475)(1174,1464)(1178,1452)
	(1182,1441)(1185,1429)(1188,1416)
	(1190,1404)(1192,1391)(1194,1377)
	(1195,1363)(1197,1348)(1199,1332)
	(1200,1316)(1201,1301)(1203,1285)
	(1205,1270)(1206,1256)(1208,1242)
	(1210,1229)(1213,1216)(1216,1202)
	(1219,1188)(1223,1174)(1228,1160)
	(1234,1145)(1239,1131)(1245,1117)
	(1251,1103)(1257,1090)(1263,1077)
	(1269,1066)(1275,1054)(1281,1042)
	(1287,1031)(1293,1018)(1300,1003)
	(1308,987)(1317,969)(1327,951)
	(1336,933)(1343,918)(1348,909)
	(1350,905)(1350,904)
\path(4946,2708)(4946,2705)(4947,2697)
	(4949,2685)(4951,2666)(4954,2644)
	(4957,2618)(4960,2590)(4964,2563)
	(4967,2536)(4970,2511)(4973,2489)
	(4976,2468)(4979,2448)(4982,2431)
	(4985,2414)(4987,2398)(4990,2383)
	(4994,2363)(4998,2344)(5002,2324)
	(5007,2304)(5012,2284)(5017,2264)
	(5022,2245)(5027,2226)(5032,2208)
	(5037,2191)(5042,2175)(5047,2161)
	(5051,2148)(5056,2136)(5062,2120)
	(5068,2106)(5074,2092)(5080,2079)
	(5087,2068)(5094,2058)(5100,2050)
	(5106,2044)(5112,2039)(5118,2036)
	(5123,2034)(5129,2033)(5134,2033)
	(5140,2034)(5145,2035)(5151,2037)
	(5157,2040)(5162,2044)(5167,2048)
	(5171,2053)(5174,2058)(5177,2064)
	(5178,2064)(5181,2073)(5184,2083)
	(5186,2095)(5187,2107)(5186,2118)
	(5184,2129)(5181,2138)(5178,2145)
	(5174,2149)(5171,2153)(5167,2155)
	(5162,2157)(5157,2159)(5152,2160)
	(5147,2160)(5142,2159)(5137,2158)
	(5132,2156)(5128,2154)(5124,2152)
	(5124,2151)(5119,2147)(5115,2140)
	(5111,2131)(5108,2120)(5106,2107)
	(5105,2091)(5105,2074)(5106,2055)
	(5106,2043)(5107,2029)(5108,2015)
	(5110,1999)(5111,1982)(5113,1964)
	(5114,1945)(5115,1926)(5117,1907)
	(5118,1888)(5118,1869)(5119,1851)
	(5119,1834)(5118,1817)(5117,1801)
	(5116,1784)(5114,1768)(5112,1751)
	(5109,1734)(5106,1717)(5103,1701)
	(5099,1685)(5095,1670)(5091,1656)
	(5087,1643)(5083,1631)(5078,1621)
	(5074,1611)(5068,1599)(5061,1588)
	(5054,1577)(5047,1568)(5040,1561)
	(5032,1554)(5025,1549)(5018,1546)
	(5011,1543)(5005,1542)(4999,1541)
	(4993,1542)(4987,1543)(4981,1545)
	(4976,1548)(4970,1552)(4966,1557)
	(4962,1562)(4958,1568)(4955,1573)
	(4952,1581)(4950,1590)(4949,1600)
	(4949,1611)(4949,1621)(4950,1630)
	(4952,1638)(4955,1645)(4959,1651)
	(4964,1656)(4969,1661)(4976,1664)
	(4983,1666)(4990,1667)(4998,1668)
	(5005,1667)(5011,1666)(5018,1664)
	(5025,1661)(5032,1658)(5039,1653)
	(5046,1648)(5052,1641)(5058,1634)
	(5063,1626)(5068,1617)(5071,1609)
	(5074,1599)(5077,1589)(5079,1577)
	(5082,1565)(5083,1551)(5084,1538)
	(5084,1524)(5084,1510)(5083,1496)
	(5082,1483)(5080,1470)(5078,1457)
	(5075,1444)(5071,1430)(5067,1415)
	(5062,1401)(5056,1386)(5050,1370)
	(5044,1356)(5037,1341)(5031,1327)
	(5024,1314)(5018,1301)(5011,1288)
	(5005,1275)(4998,1261)(4990,1247)
	(4983,1232)(4976,1217)(4969,1201)
	(4962,1186)(4956,1171)(4951,1157)
	(4947,1143)(4943,1130)(4939,1116)
	(4937,1102)(4934,1087)(4932,1072)
	(4930,1057)(4929,1042)(4928,1027)
	(4928,1014)(4927,1001)(4927,990)
	(4927,979)(4927,970)(4927,955)
	(4927,941)(4927,928)(4927,917)
	(4927,909)(4927,908)
\path(8100,2179)(8101,2179)(8104,2179)
	(8114,2179)(8130,2179)(8151,2179)
	(8178,2179)(8206,2179)(8235,2179)
	(8262,2179)(8288,2178)(8311,2178)
	(8332,2177)(8351,2176)(8368,2175)
	(8384,2174)(8400,2173)(8418,2171)
	(8435,2169)(8453,2166)(8470,2163)
	(8488,2160)(8506,2156)(8523,2152)
	(8540,2147)(8556,2143)(8571,2138)
	(8585,2133)(8598,2128)(8610,2122)
	(8621,2117)(8634,2110)(8646,2103)
	(8658,2095)(8669,2087)(8680,2077)
	(8690,2068)(8699,2058)(8707,2048)
	(8713,2038)(8718,2028)(8722,2018)
	(8725,2008)(8726,1995)(8727,1983)
	(8726,1969)(8724,1955)(8720,1940)
	(8716,1926)(8711,1913)(8705,1900)
	(8699,1889)(8694,1879)(8688,1870)
	(8681,1861)(8674,1853)(8667,1846)
	(8659,1840)(8651,1835)(8643,1831)
	(8634,1830)(8626,1829)(8618,1830)
	(8610,1832)(8603,1835)(8595,1839)
	(8587,1844)(8580,1849)(8573,1856)
	(8566,1862)(8561,1870)(8558,1877)
	(8556,1884)(8555,1891)(8556,1897)
	(8557,1897)(8559,1902)(8562,1907)
	(8565,1912)(8570,1916)(8576,1921)
	(8583,1924)(8590,1928)(8598,1931)
	(8607,1933)(8616,1934)(8624,1934)
	(8633,1934)(8642,1932)(8650,1930)
	(8658,1927)(8666,1922)(8676,1915)
	(8686,1905)(8695,1894)(8704,1880)
	(8713,1865)(8721,1849)(8728,1832)
	(8734,1815)(8738,1799)(8742,1783)
	(8744,1769)(8746,1755)(8746,1743)
	(8745,1731)(8743,1720)(8740,1709)
	(8735,1698)(8730,1688)(8724,1679)
	(8718,1671)(8710,1663)(8703,1657)
	(8695,1650)(8687,1645)(8679,1639)
	(8670,1634)(8661,1628)(8650,1623)
	(8639,1617)(8626,1612)(8613,1606)
	(8600,1601)(8586,1596)(8572,1592)
	(8558,1587)(8543,1583)(8530,1579)
	(8516,1575)(8501,1571)(8485,1567)
	(8468,1563)(8450,1559)(8432,1555)
	(8414,1551)(8396,1548)(8378,1544)
	(8361,1541)(8344,1537)(8328,1535)
	(8313,1532)(8298,1529)(8283,1527)
	(8267,1525)(8250,1522)(8232,1520)
	(8212,1517)(8190,1515)(8168,1512)
	(8146,1509)(8127,1507)(8113,1505)
	(8104,1504)(8101,1504)(8100,1504)
\put(1350,604){\makebox(0,0)[lb]{\smash{{{\SetFigFont{10}{12.0}{\rmdefault}{\mddefault}{\updefault}$X$}}}}}
\put(0,1954){\makebox(0,0)[lb]{\smash{{{\SetFigFont{10}{12.0}{\rmdefault}{\mddefault}{\updefault}$Y$}}}}}
\put(4875,2854){\makebox(0,0)[lb]{\smash{{{\SetFigFont{10}{12.0}{\rmdefault}{\mddefault}{\updefault}$Y$}}}}}
\put(7800,2104){\makebox(0,0)[lb]{\smash{{{\SetFigFont{10}{12.0}{\rmdefault}{\mddefault}{\updefault}$Y$}}}}}
\put(4800,604){\makebox(0,0)[lb]{\smash{{{\SetFigFont{10}{12.0}{\rmdefault}{\mddefault}{\updefault}$X$}}}}}
\put(7800,1429){\makebox(0,0)[lb]{\smash{{{\SetFigFont{10}{12.0}{\rmdefault}{\mddefault}{\updefault}$X$}}}}}
\put(450,1579){\makebox(0,0)[lb]{\smash{{{\SetFigFont{10}{12.0}{\rmdefault}{\mddefault}{\updefault}$g$}}}}}
\put(1200,1654){\makebox(0,0)[lb]{\smash{{{\SetFigFont{10}{12.0}{\rmdefault}{\mddefault}{\updefault}$f$}}}}}
\put(4350,1804){\makebox(0,0)[lb]{\smash{{{\SetFigFont{10}{12.0}{\rmdefault}{\mddefault}{\updefault}$g$}}}}}
\put(8175,1804){\makebox(0,0)[lb]{\smash{{{\SetFigFont{10}{12.0}{\rmdefault}{\mddefault}{\updefault}$g$}}}}}
\put(5325,1729){\makebox(0,0)[lb]{\smash{{{\SetFigFont{10}{12.0}{\rmdefault}{\mddefault}{\updefault}$f$}}}}}
\put(8925,1804){\makebox(0,0)[lb]{\smash{{{\SetFigFont{10}{12.0}{\rmdefault}{\mddefault}{\updefault}$f$}}}}}
\put(4200,3379){\makebox(0,0)[lb]{\smash{{{\SetFigFont{10}{12.0}{\rmdefault}{\mddefault}{\updefault}$q$}}}}}
\put(8100,3379){\makebox(0,0)[lb]{\smash{{{\SetFigFont{10}{12.0}{\rmdefault}{\mddefault}{\updefault}$q$}}}}}
\put(6675,679){\makebox(0,0)[lb]{\smash{{{\SetFigFont{10}{12.0}{\rmdefault}{\mddefault}{\updefault}$p$}}}}}
\put(10575,679){\makebox(0,0)[lb]{\smash{{{\SetFigFont{10}{12.0}{\rmdefault}{\mddefault}{\updefault}$p$}}}}}
\put(1309,82){\makebox(0,0)[lb]{\smash{{{\SetFigFont{10}{8.0}{\rmdefault}{\mddefault}{\updefault}a}}}}}
\put(5209,83){\makebox(0,0)[lb]{\smash{{{\SetFigFont{10}{8.0}{\rmdefault}{\mddefault}{\updefault}b}}}}}
\put(9090,82){\makebox(0,0)[lb]{\smash{{{\SetFigFont{10}{8.0}{\rmdefault}{\mddefault}{\updefault}c}}}}}
\put(2884,1620){\makebox(0,0)[lb]{\smash{{{\SetFigFont{10}{12.0}{\rmdefault}{\mddefault}{\updefault}$\phi$}}}}}
\put(3053,570){\makebox(0,0)[lb]{\smash{{{\SetFigFont{10}{12.0}{\rmdefault}{\mddefault}{\updefault}$\rho$}}}}}
\put(189,507){\makebox(0,0)[lb]{\smash{{{\SetFigFont{10}{12.0}{\rmdefault}{\mddefault}{\updefault}$O$}}}}}
\put(4112,507){\makebox(0,0)[lb]{\smash{{{\SetFigFont{10}{12.0}{\rmdefault}{\mddefault}{\updefault}$O$}}}}}
\put(8012,507){\makebox(0,0)[lb]{\smash{{{\SetFigFont{10}{12.0}{\rmdefault}{\mddefault}{\updefault}$O$}}}}}
\end{picture}
}
$$
\caption{}
\end{figure}
{\bf Case 1.} Choose $g$ so that it contains a circular arc
whose endpoints lie on the linear segments passing through $X$ and $Y$
(see figure 2a).
 We work in polar
coordinates  $(\rho,\phi)$ with the origin at $O$. Let
$(\rho(t),\phi(t)),\; t \in [0;\pi/2]$ be a (re)parameterization
of $f$.
Consider an arbitrarily small perturbation of $g$ of the form
$h = (u(t),t), \; t \in [0;\pi/2]$.
Then $u(\phi(t)) \leq \rho (t)$
for all $t$. Consider the homotopy
$$\theta (s,t) = ((1-s)\rho(t) + su(\phi(t)),\phi(t)), s \in [0;1]$$
between $f$ and the curve $h'(t) = (u(\phi(t)),\phi(t))$.
In polar coordinates the area form $\omega$ is given by
$\rho d\rho \wedge d\phi$ so
$$|\theta^*\omega| =
((1-s)\rho(t) + su(\phi(t)) |u(\phi(t)) - \rho(t)|\; 
|{\dot \phi}(t)| \; ds \wedge dt
\leq \rho(t) \sqrt{2} |{\dot \phi}(t)| ds \wedge dt.$$
Hence the area of homotopy $\theta$  can be estimated as follows:
$$\int |\theta^*\omega| \leq \sqrt{2} \int_0^{\pi/2} \rho 
|{\dot \phi}| dt
\leq \sqrt{2} \Le (f).$$ 
Note that $h$ and $h'$ are homotopic 
with fixed end points in $\R^2 \setminus \Z^2$ by a homotopy of area 
$0$. Combining it with $\theta$, we get the desired homotopy.

Further,
$$\Le (f) \geq \sqrt{OX^2 + OY^2} \geq
(OX + OY)/\sqrt{2} $$
$$\geq
(OX +OY - 2\epsilon + \pi\epsilon/2)/\sqrt{2} = \Le (g)/\sqrt{2}.$$

\medskip
\noindent
{\bf Case 2.} Choose $g$ so that $\Le (g) \leq 2$.
We work in Euclidean coordinates $(p,q)$ with the origin at $O$.
After  a (re)parameterization  we write
$f(t) = (p(t),q(t)), \; t \in [0;1]$.
Perturb $g$ to a curve $h(t) = (u(t),t),\; t \in [0;1]$.
Consider the homotopy
$$\theta (s,t) = ((1-s)p(t)+su(q(t)),q(t)), s \in [0;1]$$
between $f$ and the curve
$h'(t) = (u(q(t)),q(t))$.
Since $\omega = dp \wedge dq$
we
have
$$|\theta^*\omega| = |u(q(t))-p(t)|\;|{\dot q}(t)|ds\wedge dt
\leq |{\dot q}(t)| ds \wedge dt,$$
so the area of $\theta$ does not exceed
$$\int_0^1 |{\dot q}(t)| dt \leq \Le (f).$$
Note that $h$ and $h'$ are homotopic 
with fixed end points in $\R^2 \setminus \Z^2$ by a homotopy of area 
$0$. Combining it with $\theta$, we get the desired homotopy.

Further,
$\Le (f) \geq 1 \geq \Le (g)/2$.

\medskip
\noindent
{\bf Case 3.} Choose $g = XY$. In Euclidean coordinates
$(p,q)$ with the origin at $O$
 one has $g(t) = (0,t)$ and $f(t) = (p(t),q(t))$,
where $t \in [q_0;q_1]$. Here $X = (0,q_0)$ and $Y = (0,q_1)$.
One computes that the area of the
homotopy $(sp(t),q(t)), \; s \in [0;1]$ does not exceed $\Le (f)$.
One can combine it with the obvious homotopy with fixed
end points in $\R^2 \setminus \Z^2$ of area $0$ to get the desired
homotopy between $f$ and $g$. Further,
$\Le (f) \geq \Le (g)$.

This finishes off the proof of the Lemma  in the case when
$\alpha \cap G_{\epsilon} \neq \emptyset$. Suppose now that
this intersection is empty. Then $\alpha$ is entirely contained
in some fundamental domain $Q$ of $\R^2 \setminus B(\epsilon)$.
Arguing exactly as in Step 3 above we get the desired homotopy.
This completes the proof.
$\QED$.


 \section{Combinatorial isoperimetric inequality}

Let $w \in \F_3$ be a cyclically reduced word. A {\it simple conjugate}
of $w$ is a word of the form $\beta\alpha$ where
$\alpha$ is an initial sub-word of $w$ and
$w = \alpha\beta$.
Denote by $R$ the set of 10 words consisting of simple conjugates of
$r^{\pm 1}$, that is
$R = \{r_1^{\pm 1}, ... , r_5^{\pm 1}\}$, where
$r_1 = r = aba^{-1}b^{-1}c$, $r_2 = ba^{-1}b^{-1}ca$, $r_3 = a^{-1}b^{-1}cab$,
$r_4 = b^{-1}caba^{-1}$ and $r_5 = caba^{-1}b^{-1}$.
A {\it piece} is an initial sub-word of some $r \in R$.

\medskip

\begin{thm}
\prlabel{piece}
Every cyclically reduced word $w \in N$
contains a piece
of length 4.
\end{thm}

Note that the small cancellation theory guarantees
existence of such a piece of length 3 only.

\medskip
\noindent
{\bf Proof of Theorem \ref{ineq2}:} We prove $A(w) \leq l(w)/2$
by the induction in the word length $|w|$ of $w \in N$.
If $|w| = 5$ then $w \in R$ (one can see this directly
on the Cayley complex, see figure 3 below). Hence
$A(w) = 1$ while $l(w) = 4$ so the inequality
holds.  Take $w \in N$ and assume that our inequality
is already proved for all words
of length $< |w|$. By Theorem \ref{piece} $w$ contains a piece, say $u$,
of length 4. The piece $u$ is the initial
sub-word of some unique $s \in R$.
Define $u_*$ by
$s = uu_*^{-1}$. Consider a new word $w' \in N$ which is obtained from $w$
by replacing $u$ with $u_*$ and cyclic reduction. Clearly, one has
that $A(w) \leq A(w') +1$, $l(w') \leq l(w) -2$ and $|w'| \leq |w| -3$.
By the inductive assumption, $A(w') \leq l(w')/2$. Combining the
inequalities we get $A(w) \leq l(w)/2$ as required.
$\QED$
\begin{figure}[h]
\label{fig3}
$$
\setlength{\unitlength}{0.00069167in}
\begingroup\makeatletter\ifx\SetFigFont\undefined%
\gdef\SetFigFont#1#2#3#4#5{%
  \reset@font\fontsize{#1}{#2pt}%
  \fontfamily{#3}\fontseries{#4}\fontshape{#5}%
  \selectfont}%
\fi\endgroup%
{\renewcommand{\dashlinestretch}{30}
\begin{picture}(5855,6069)(0,-10)
\put(4668.401,790.833){\arc{2162.683}{4.1524}{5.3464}}
\path(4184.675,1791.131)(4094.000,1707.000)(4213.599,1738.562)
\put(5836.763,2950.816){\arc{2242.131}{2.0012}{3.2634}}
\path(4746.066,2965.291)(4724.000,3087.000)(4686.195,2969.226)
\put(981,4242){\ellipse{376}{376}}
\put(2406,5217){\ellipse{376}{376}}
\put(3906,4242){\ellipse{376}{376}}
\put(3019,1397){\ellipse{376}{376}}
\put(1144,2222){\ellipse{376}{376}}
\path(269,5142)(44,4092)
\path(44,4092)(569,3342)
\path(569,3342)(2444,3342)
\path(1244,5142)(2444,3342)
\path(2594,3492)(1469,5217)
\path(1469,5217)(2069,5967)
\path(2069,5967)(2894,5967)
\path(2894,5967)(3269,5367)
\path(3269,5367)(2594,3492)
\path(2744,3342)(3419,5292)
\path(3419,5292)(4769,5292)
\path(5219,4317)(4919,3342)
\path(4919,3342)(2744,3342)
\path(2444,3192)(569,3192)
\path(569,3192)(44,1992)
\path(44,1992)(419,1167)
\path(419,1167)(1769,1092)
\path(1769,1092)(2444,3192)
\path(2669,3192)(1994,1017)
\path(2744,42)(3719,117)
\path(3719,117)(4244,1242)
\path(4244,1242)(2669,3192)
\blacken\path(58,4495)(158,4630)(206,4470)(58,4495)
\path(58,4495)(158,4630)(206,4470)(58,4495)
\path(1244,5142)(269,5142)
\blacken\path(776,5217)(925,5142)(775,5067)(776,5217)
\path(776,5217)(925,5142)(775,5067)(776,5217)
\blacken\path(1806,4424)(1844,4261)(1689,4328)(1806,4424)
\path(1806,4424)(1844,4261)(1689,4328)(1806,4424)
\blacken\path(2445,6042)(2594,5967)(2444,5892)(2445,6042)
\path(2445,6042)(2594,5967)(2444,5892)(2445,6042)
\blacken\path(4020,5367)(4169,5292)(4019,5217)(4020,5367)
\path(4020,5367)(4169,5292)(4019,5217)(4020,5367)
\blacken\path(1286,1029)(1150,1130)(1310,1179)(1286,1029)
\path(1286,1029)(1150,1130)(1310,1179)(1286,1029)
\blacken\path(1634,5563)(1788,5630)(1749,5467)(1634,5563)
\path(1634,5563)(1788,5630)(1749,5467)(1634,5563)
\blacken\path(185,2537)(325,2630)(315,2463)(185,2537)
\path(185,2537)(325,2630)(315,2463)(185,2537)
\blacken\path(175,1469)(194,1636)(316,1521)(175,1469)
\path(175,1469)(194,1636)(316,1521)(175,1469)
\blacken\path(2069,4131)(2088,4298)(2210,4183)(2069,4131)
\path(2069,4131)(2088,4298)(2210,4183)(2069,4131)
\blacken\path(3081,5793)(3119,5630)(2964,5697)(3081,5793)
\path(3081,5793)(3119,5630)(2964,5697)(3081,5793)
\path(4769,5292)(5219,4317)
\blacken\path(5022,4897)(5031,4730)(4890,4823)(5022,4897)
\path(5022,4897)(5031,4730)(4890,4823)(5022,4897)
\blacken\path(3092,4656)(2969,4542)(2949,4710)(3092,4656)
\path(3092,4656)(2969,4542)(2949,4710)(3092,4656)
\blacken\path(2922,4126)(3044,4242)(3063,4076)(2922,4126)
\path(2922,4126)(3044,4242)(3063,4076)(2922,4126)
\blacken\path(2247,2101)(2369,2217)(2388,2051)(2247,2101)
\path(2247,2101)(2369,2217)(2388,2051)(2247,2101)
\blacken\path(5211,3962)(5088,3848)(5068,4016)(5211,3962)
\path(5211,3962)(5088,3848)(5068,4016)(5211,3962)
\blacken\path(2098,1862)(1975,1748)(1955,1916)(2098,1862)
\path(2098,1862)(1975,1748)(1955,1916)(2098,1862)
\blacken\path(3466,2348)(3475,2180)(3334,2273)(3466,2348)
\path(3466,2348)(3475,2180)(3334,2273)(3466,2348)
\blacken\path(3929,3259)(3794,3361)(3956,3410)(3929,3259)
\path(3929,3259)(3794,3361)(3956,3410)(3929,3259)
\blacken\path(307,3590)(269,3755)(425,3689)(307,3590)
\path(307,3590)(269,3755)(425,3689)(307,3590)
\blacken\path(4104,753)(3963,661)(3972,830)(4104,753)
\path(4104,753)(3963,661)(3972,830)(4104,753)
\blacken\path(3336,12)(3175,61)(3310,163)(3336,12)
\path(3336,12)(3175,61)(3310,163)(3336,12)
\Thicklines
\path(269,5142)(44,4092)
\path(269,5142)(44,4092)
\path(569,3192)(44,1992)
\path(569,3192)(44,1992)
\path(2669,3192)(4244,1242)
\path(2669,3192)(4244,1242)
\path(2744,3342)(4919,3342)
\path(2744,3342)(4919,3342)
\path(1469,5217)(2069,5967)
\path(1469,5217)(2069,5967)
\thinlines
\blacken\path(1818,3259)(1667,3336)(1818,3412)(1818,3259)
\path(1818,3259)(1667,3336)(1818,3412)(1818,3259)
\blacken\path(1411,3283)(1560,3208)(1410,3133)(1411,3283)
\path(1411,3283)(1560,3208)(1410,3133)(1411,3283)
\path(3194,3042)(2894,3192)
\path(2744,42)(1994,1017)
\blacken\path(2303,474)(2294,642)(2435,549)(2303,474)
\path(2303,474)(2294,642)(2435,549)(2303,474)
\blacken\path(2911,3111)(2819,3252)(2987,3242)(2911,3111)
\path(2911,3111)(2819,3252)(2987,3242)(2911,3111)
\path(4124,1692)(4169,1827)
\path(4724,3072)(4634,2952)
\put(869,4827){\makebox(0,0)[lb]{\smash{{{\SetFigFont{12}{14.4}{\rmdefault}{\mddefault}{\updefault}$a$}}}}}
\put(2429,5712){\makebox(0,0)[lb]{\smash{{{\SetFigFont{12}{14.4}{\rmdefault}{\mddefault}{\updefault}$a$}}}}}
\put(3239,4287){\makebox(0,0)[lb]{\smash{{{\SetFigFont{12}{14.4}{\rmdefault}{\mddefault}{\updefault}$a$}}}}}
\put(3749,717){\makebox(0,0)[lb]{\smash{{{\SetFigFont{12}{14.4}{\rmdefault}{\mddefault}{\updefault}$a$}}}}}
\put(1589,2997){\makebox(0,0)[lb]{\smash{{{\SetFigFont{12}{14.4}{\rmdefault}{\mddefault}{\updefault}$a$}}}}}
\put(1544,4212){\makebox(0,0)[lb]{\smash{{{\SetFigFont{12}{14.4}{\rmdefault}{\mddefault}{\updefault}$b$}}}}}
\put(2894,5472){\makebox(0,0)[lb]{\smash{{{\SetFigFont{12}{14.4}{\rmdefault}{\mddefault}{\updefault}$b$}}}}}
\put(3929,5022){\makebox(0,0)[lb]{\smash{{{\SetFigFont{12}{14.4}{\rmdefault}{\mddefault}{\updefault}$b$}}}}}
\put(3194,282){\makebox(0,0)[lb]{\smash{{{\SetFigFont{12}{14.4}{\rmdefault}{\mddefault}{\updefault}$b$}}}}}
\put(1979,2187){\makebox(0,0)[lb]{\smash{{{\SetFigFont{12}{14.4}{\rmdefault}{\mddefault}{\updefault}$b$}}}}}
\put(1424,3447){\makebox(0,0)[lb]{\smash{{{\SetFigFont{12}{14.4}{\rmdefault}{\mddefault}{\updefault}$a^{-1}$}}}}}
\put(2644,4782){\makebox(0,0)[lb]{\smash{{{\SetFigFont{12}{14.4}{\rmdefault}{\mddefault}{\updefault}$a^{-1}$}}}}}
\put(4669,4542){\makebox(0,0)[lb]{\smash{{{\SetFigFont{12}{14.4}{\rmdefault}{\mddefault}{\updefault}$a^{-1}$}}}}}
\put(2459,537){\makebox(0,0)[lb]{\smash{{{\SetFigFont{12}{14.4}{\rmdefault}{\mddefault}{\updefault}$a^{-1}$}}}}}
\put(1004,1212){\makebox(0,0)[lb]{\smash{{{\SetFigFont{12}{14.4}{\rmdefault}{\mddefault}{\updefault}$a^{-1}$}}}}}
\put(404,3762){\makebox(0,0)[lb]{\smash{{{\SetFigFont{12}{14.4}{\rmdefault}{\mddefault}{\updefault}$b^{-1}$}}}}}
\put(1994,4542){\makebox(0,0)[lb]{\smash{{{\SetFigFont{12}{14.4}{\rmdefault}{\mddefault}{\updefault}$b^{-1}$}}}}}
\put(4669,3732){\makebox(0,0)[lb]{\smash{{{\SetFigFont{12}{14.4}{\rmdefault}{\mddefault}{\updefault}$b^{-1}$}}}}}
\put(2489,2202){\makebox(0,0)[lb]{\smash{{{\SetFigFont{12}{14.4}{\rmdefault}{\mddefault}{\updefault}$b^{-1}$}}}}}
\put(239,1677){\makebox(0,0)[lb]{\smash{{{\SetFigFont{12}{14.4}{\rmdefault}{\mddefault}{\updefault}$b^{-1}$}}}}}
\put(299,4542){\makebox(0,0)[lb]{\smash{{{\SetFigFont{12}{14.4}{\rmdefault}{\mddefault}{\updefault}$c$}}}}}
\put(1904,5517){\makebox(0,0)[lb]{\smash{{{\SetFigFont{12}{14.4}{\rmdefault}{\mddefault}{\updefault}$c$}}}}}
\put(4049,3432){\makebox(0,0)[lb]{\smash{{{\SetFigFont{12}{14.4}{\rmdefault}{\mddefault}{\updefault}$c$}}}}}
\put(3344,1992){\makebox(0,0)[lb]{\smash{{{\SetFigFont{12}{14.4}{\rmdefault}{\mddefault}{\updefault}$c$}}}}}
\put(479,2562){\makebox(0,0)[lb]{\smash{{{\SetFigFont{12}{14.4}{\rmdefault}{\mddefault}{\updefault}$c$}}}}}
\put(5579,1662){\makebox(0,0)[lb]{\smash{{{\SetFigFont{12}{14.4}{\rmdefault}{\mddefault}{\updefault}$\partial
K$}}}}}
\put(854,4167){\makebox(0,0)[lb]{\smash{{{\SetFigFont{12}{14.4}{\rmdefault}{\mddefault}{\updefault}$r_3$}}}}}
\put(1004,2127){\makebox(0,0)[lb]{\smash{{{\SetFigFont{12}{14.4}{\rmdefault}{\mddefault}{\updefault}$r_2$}}}}}
\put(2300,5127){\makebox(0,0)[lb]{\smash{{{\SetFigFont{12}{14.4}{\rmdefault}{\mddefault}{\updefault}$r_4$}}}}}
\put(3794,4182){\makebox(0,0)[lb]{\smash{{{\SetFigFont{12}{14.4}{\rmdefault}{\mddefault}{\updefault}$r_1$}}}}}
\put(2909,1317){\makebox(0,0)[lb]{\smash{{{\SetFigFont{12}{14.4}{\rmdefault}{\mddefault}{\updefault}$r_5$}}}}}
\put(3381,2892){\makebox(0,0)[lb]{\smash{{{\SetFigFont{10}{12.0}{\rmdefault}{\mddefault}{\updefault}$O$}}}}}
\end{picture}
}
$$
\caption{}
\end{figure}

\medskip
\noindent
It remains to prove Theorem \ref{piece}. Denote by $K$ the
Cayley complex of the group $\F_2 = \F_3 /N$ associated to the
presentation $\F_2 = <a,b,c|aba^{-1}b^{-1}c>$. It is easy to see
that $K$ is a 2-dimensional  simply connected PL-manifold with boundary
(in fact $K$ is homeomorphic to the universal cover of $\T^2 \setminus D^2$).
Every vertex, say $O$ of $K$ lies on the boundary of $K$.
The union of all faces which have $O$ as a vertex
is presented on figure 3.

\medskip
\noindent
{\bf Proof of Theorem \ref{piece} :} Let $w$ be a cyclically reduced
word representing an element from $N$ . Consider the
cycle on $K$ which corresponds to $w$
and starts at 1. This cycle has an embedded subcycle, say $\gamma$,
which corresponds to some sub-word $w'$ of $w$
(see figure 4).
\begin{figure}[h]\label{fig4}
$$
\setlength{\unitlength}{0.00073333in}
\begingroup\makeatletter\ifx\SetFigFont\undefined%
\gdef\SetFigFont#1#2#3#4#5{%
  \reset@font\fontsize{#1}{#2pt}%
  \fontfamily{#3}\fontseries{#4}\fontshape{#5}%
  \selectfont}%
\fi\endgroup%
{\renewcommand{\dashlinestretch}{30}
\begin{picture}(4128,1958)(0,-10)
\path(519,788)(744,806)(669,619)
\path(1081,1481)(1119,1706)(1250,1556)
\path(2769,1669)(2881,1481)(2675,1500)
\path(1531,338)(1531,544)(1700,450)
\path(12,94)(15,97)(21,103)
	(32,114)(46,128)(64,146)
	(83,165)(103,184)(122,203)
	(139,221)(155,237)(170,252)
	(184,265)(197,278)(210,291)
	(222,303)(233,314)(244,325)
	(256,337)(268,349)(280,361)
	(293,374)(306,387)(320,401)
	(333,414)(347,428)(361,442)
	(375,456)(388,469)(401,483)
	(415,496)(428,509)(440,521)
	(453,534)(466,547)(480,561)
	(493,574)(508,588)(522,603)
	(538,618)(554,633)(570,649)
	(586,665)(603,680)(619,696)
	(636,711)(652,725)(668,739)
	(683,753)(698,766)(713,778)
	(728,790)(743,802)(758,814)
	(774,826)(789,838)(806,850)
	(822,862)(839,874)(856,886)
	(873,898)(890,910)(907,922)
	(923,934)(939,945)(954,956)
	(969,966)(983,977)(996,987)
	(1009,997)(1022,1007)(1035,1017)
	(1048,1028)(1061,1039)(1074,1050)
	(1088,1062)(1102,1074)(1116,1086)
	(1131,1098)(1145,1111)(1160,1123)
	(1175,1136)(1190,1148)(1205,1160)
	(1220,1172)(1235,1183)(1250,1195)
	(1266,1206)(1280,1217)(1295,1227)
	(1311,1238)(1327,1249)(1344,1260)
	(1362,1272)(1381,1284)(1400,1296)
	(1419,1308)(1439,1320)(1459,1332)
	(1478,1343)(1498,1355)(1517,1366)
	(1535,1376)(1553,1386)(1570,1396)
	(1587,1405)(1603,1414)(1619,1422)
	(1636,1431)(1653,1440)(1670,1448)
	(1687,1456)(1704,1465)(1721,1473)
	(1739,1481)(1756,1489)(1774,1496)
	(1792,1503)(1809,1510)(1827,1517)
	(1844,1523)(1860,1529)(1877,1535)
	(1893,1540)(1909,1545)(1925,1550)
	(1941,1555)(1958,1559)(1976,1564)
	(1994,1568)(2013,1573)(2033,1578)
	(2053,1582)(2074,1586)(2095,1590)
	(2116,1594)(2137,1598)(2158,1602)
	(2178,1605)(2198,1607)(2217,1610)
	(2236,1612)(2254,1614)(2272,1615)
	(2292,1617)(2311,1618)(2332,1619)
	(2352,1620)(2373,1620)(2394,1620)
	(2416,1619)(2437,1619)(2458,1617)
	(2478,1616)(2498,1615)(2517,1613)
	(2535,1611)(2552,1608)(2569,1606)
	(2584,1603)(2600,1601)(2615,1598)
	(2631,1594)(2647,1591)(2663,1586)
	(2679,1582)(2696,1577)(2713,1571)
	(2729,1565)(2746,1559)(2762,1552)
	(2778,1546)(2794,1539)(2810,1531)
	(2825,1524)(2841,1516)(2856,1508)
	(2873,1499)(2890,1489)(2907,1479)
	(2925,1468)(2944,1457)(2963,1445)
	(2982,1432)(3001,1420)(3019,1407)
	(3036,1395)(3053,1383)(3069,1371)
	(3083,1359)(3097,1348)(3110,1337)
	(3123,1325)(3136,1313)(3149,1300)
	(3161,1287)(3172,1273)(3183,1259)
	(3194,1245)(3203,1230)(3212,1215)
	(3220,1200)(3227,1186)(3233,1171)
	(3239,1156)(3244,1141)(3248,1127)
	(3252,1113)(3256,1098)(3259,1082)
	(3263,1065)(3266,1048)(3269,1030)
	(3272,1011)(3275,993)(3278,974)
	(3280,956)(3282,938)(3284,921)
	(3285,904)(3287,888)(3288,872)
	(3289,857)(3290,841)(3290,825)
	(3291,809)(3291,792)(3291,775)
	(3291,757)(3291,739)(3290,721)
	(3289,703)(3287,685)(3285,667)
	(3283,650)(3281,632)(3278,615)
	(3275,597)(3272,581)(3269,565)
	(3265,548)(3260,530)(3255,511)
	(3250,492)(3244,473)(3237,453)
	(3229,433)(3221,414)(3213,394)
	(3204,375)(3194,357)(3184,340)
	(3174,323)(3163,307)(3152,292)
	(3141,278)(3128,265)(3115,252)
	(3101,239)(3086,226)(3070,213)
	(3054,201)(3036,189)(3018,178)
	(2999,166)(2980,156)(2961,146)
	(2942,136)(2924,127)(2905,119)
	(2888,111)(2871,104)(2854,97)
	(2838,91)(2819,84)(2801,77)
	(2783,71)(2764,65)(2745,58)
	(2725,53)(2706,47)(2686,42)
	(2666,37)(2647,33)(2628,29)
	(2610,25)(2592,22)(2574,20)
	(2557,17)(2541,16)(2524,14)
	(2507,13)(2489,12)(2471,12)
	(2452,12)(2433,12)(2413,12)
	(2392,13)(2372,15)(2351,16)
	(2331,18)(2311,20)(2291,23)
	(2272,26)(2253,28)(2235,31)
	(2218,34)(2200,38)(2183,41)
	(2165,45)(2146,49)(2127,54)
	(2108,59)(2088,64)(2069,70)
	(2049,76)(2030,82)(2012,88)
	(1994,95)(1978,101)(1962,108)
	(1946,115)(1932,121)(1919,128)
	(1904,137)(1889,145)(1875,155)
	(1861,165)(1847,175)(1833,187)
	(1819,199)(1805,211)(1792,224)
	(1778,237)(1765,251)(1752,264)
	(1740,278)(1728,291)(1715,305)
	(1703,319)(1692,331)(1681,344)
	(1669,358)(1656,372)(1643,387)
	(1630,403)(1616,419)(1603,436)
	(1589,453)(1575,470)(1562,488)
	(1548,505)(1536,522)(1524,539)
	(1512,555)(1501,571)(1491,587)
	(1481,603)(1472,619)(1463,635)
	(1454,652)(1445,668)(1436,686)
	(1428,704)(1420,722)(1411,741)
	(1403,759)(1396,778)(1388,797)
	(1381,815)(1374,833)(1368,851)
	(1361,867)(1355,884)(1350,900)
	(1344,916)(1338,932)(1332,948)
	(1326,964)(1320,981)(1314,999)
	(1307,1018)(1300,1037)(1294,1057)
	(1287,1078)(1280,1100)(1272,1122)
	(1265,1144)(1258,1167)(1252,1190)
	(1245,1214)(1238,1238)(1231,1262)
	(1225,1288)(1220,1307)(1215,1328)
	(1209,1349)(1204,1372)(1198,1397)
	(1192,1424)(1185,1452)(1178,1484)
	(1171,1517)(1163,1553)(1154,1592)
	(1145,1633)(1136,1675)(1127,1717)
	(1118,1759)(1109,1798)(1102,1834)
	(1095,1865)(1090,1890)(1086,1909)
	(1083,1921)(1082,1928)(1081,1931)
\put(1191,1101){\makebox(0,0)[lb]{\smash{{{\SetFigFont{12}{14.4}{\rmdefault}{\mddefault}{\updefault}$\bullet$}}}}}
\put(1569,938){\makebox(0,0)[lb]{\smash{{{\SetFigFont{12}{14.4}{\rmdefault}{\mddefault}{\updefault}$P$}}}}}
\put(3575,356){\makebox(0,0)[lb]{\smash{{{\SetFigFont{12}{14.4}{\rmdefault}{\mddefault}{\updefault}$\gamma$}}}}}
\end{picture}
}
$$
\caption{}
\end{figure}
 Let $U \subset K$ be
a subcomplex (which is topologically a disc) bounded by $\gamma$.
Write $U = \Delta_1 \cup ... \cup \Delta_F$, where $\Delta_i$
are faces of $K$.

\medskip

\begin{lemma}
\prlabel{fund}
 Either $U$ consists of exactly 1 face, or there
exist at least 2 distinct faces $\Delta_i,\Delta_j$ such that
each of them intersects $\partial U$ along 4 edges.
\end{lemma}

\medskip
\noindent
Assume the Lemma. If $U$ is a face then $w' \in R$ so $w$
has a piece of length 4. Otherwise, let $\Delta$ be a face
whose intersection with $\partial U$ contains 4 edges as provided
by the Lemma. Since
$\gamma$ is embedded we conclude that these 4 edges must be
{\it consecutive} edges for $\gamma$. Thus $\gamma$ contains two
subsegments with disjoint interiors corresponding to pieces of
length 4. At most one of them contains the point $P$ (see figure 4)
as an interior point. Hence the second segment supplies us with
a piece of length 4 for $w'$, and hence for $w$. This completes
the proof of Theorem \ref{piece} modulo Lemma \ref{fund}.
$\QED$

\medskip
\noindent
{\bf Proof of Lemma \ref{fund}:} Let $V,E$ and $F$ be the number of vertices,
edges and faces of the subcomplex $U$. Denote by $E'$ the number of boundary
edges. Then $V = E'$ since every vertex of $U$ lies on the boundary of $K$
and hence on the boundary of $U$. Since every face has precisely 5 edges one has
$5F = 2E - E'$. Substituting this into the Euler formula $V - E + F = 1$
we get that $E' = 3F + 2$. But the total number of faces equals $F$, so
either there is a face with 5 boundary edges (which means $F=1$) or
there are 2 faces with at least 4 boundary edges.
This completes the proof.
$\QED$

\medskip
\noindent
{\bf 3.3 Example.} Consider a sequence of elements
$u_k \in N, \; k \in \N$ defined as
$$u_k = (b^{-1}c)^kab^ka^{-1}.$$
We claim that the ratio $A(u_k)/l(u_k)$ goes to $1/2$
as $k \to \infty$, and therefore the isoperimetric constant
in Theorem \ref{ineq2} cannot be improved.

First of let us check that $A(u_k) \leq k$. We use induction in $k$.
The word $u_1$ is contained in $R$, so its combinatorial area
equals 1. Further, $u_kr = (b^{-1}c)^{k}ab^{k+1}a^{-1}b^{-1}c$
and hence is conjugate to $u_{k+1}$. Therefore,
$A(u_{k+1}) \leq A(u_k) +1 \leq k+1$, as required.

Further, consider a homomorphism $\phi :\F_3 \to \R$
such that $\phi(a) = \phi(b) = 0$ and $\phi(c) = 1$.
Since $\phi(r) = 1$ we have that $|\phi(w)| \leq A(w)$
for every $w \in N$. Note that $\phi(u_k)=k$. Hence
$A(u_k) \geq k$.

We conclude that $A(u_k) = k$, while $l(u_k) = 2k+1$.
The claim follows.

\section{An application to dynamics}

In this section we prove Theorem \ref{growth}
by combining a method developed in \cite{P2} with
isoperimetric inequality \ref{ineq1}.

The famous Arnold conjecture proved in \cite{CZ} states
that $\psi$ has at least two fixed points, say $x$ and $y$,
with distinct projections to $\T^2$. Join them by any curve
$\gamma$ and consider a loop $\alpha: S^1 \to \R^2$ formed
by $\gamma$ and its image under $\psi$ taken with the opposite
orientation. Extend $\alpha$ to any map $f:D^2 \to \R^2$.
Its symplectic area
$$\delta(x,y,\psi) = \int_{D^2} f^*\omega$$
does not depend on the particular choice of $\gamma$ and $f$
and is called {\it the action difference} of $x$ and $y$.
Moreover, this quantity behaves nicely under iterations:
\begin{equation}
\prlabel{act}
\delta(x,y,\psi^n) = n\delta(x,y,\psi)
\end{equation}
for all $n \in \N$.
Modern developments in symplectic topology led to the following
refinement of the Conley-Zehnder result cited above: $\psi$ {\it has
a pair of fixed points with positive action difference} (see \cite{Sch}).

So take fixed points $x,y$ with $\delta(x,y,\psi)> 0$ and a
map $f$ as above. There exists a vector $v \in \Z^2$ and a
sufficiently large positive integer $m$ so that the lattice
$L = x+v+m\Z^2$ is disjoint from the image of $f$.
Then $\psi (\gamma)$ and $\gamma$ are homotopic with fixed end points
in $\R^2 \setminus L$. Since $L$ consists of fixed points of $\psi$
we see that $\psi^k(\gamma)$ is homotopic to $\psi^{k-1} (\gamma)$
with fixed end points for all $k \in \N$. Thus the loop $\alpha_n$
formed by $\gamma$ and $\psi^n (\gamma)$ is contractible
in $\R^2 \setminus L$.
Therefore Theorem \ref{ineq1} yields
\begin{equation}
\prlabel{ap}
\Ar(\alpha_n) \leq \mu \Le (\alpha_n),
\end{equation}
where the constant $\mu$ depends only on the lattice $L$ and not
on $n$.

Note now that the symplectic area of any loop does not exceed
its Euclidean area, so $$\delta(x,y,\psi^n) \leq \Ar (\alpha_n).$$
Combining this with inequalities (\ref{act}) and (\ref{ap})
we get that
$$\Le (\psi^n(\gamma)) = \Le (\alpha_n) - \Le (\gamma)$$
$$\geq \mu^{-1} \Ar (\alpha_n) - \Le (\gamma) \geq \mu^{-1}\delta(x,y,\psi^n)
-\Le (\gamma)$$
$$= n\mu^{-1}\delta(x,y,\psi) - \Le (\gamma).$$
On the other hand
$$\Le (\psi^n(\gamma)) \leq ||d\psi^n||\Le (\gamma).$$
Therefore
$||d\psi^n|| \geq \kappa n$ with
$$\kappa = (1 +\Le (\gamma))^{-1}\mu^{-1}\delta(x,y,\psi).$$
This completes the proof of Theorem \ref{growth}.

\medskip
\noindent
{\bf Remark 3.1.} Exactly the same method enables us to prove that
$||d\psi^n|| \geq \text{const}\; n$ for any {\it non-Hamiltonian}
$\Z^2$-periodic area-preserving  diffeomorphism
$\psi $ which has a fixed point.
If $x$ is such a fixed point, then one can choose $v \in \Z^2$
so that $\delta(x,x+v,\psi) > 0$. Existence of $v$ follows
from the theory of Flux homomorphism, see \cite{MS},\cite{P1}.
The rest of the argument goes through without changes.
Combining this with Theorem \ref{growth} we see that
$||d\psi^n|| \geq \text{const}\; n$
{\it for every area-preserving diffeomorphism
 $\psi \neq \text{id}$ of $\T^2$ which is isotopic to identity
and has a fixed point.}

\medskip
\noindent
{\bf Acknowledgment.} We thank Misha Kapovich,
Misha Katz, Bruce Kleiner and Joshua Maher for interesting comments.
A part of this  work was written
when both authors attended the symplectic workshop
at CRM (Montreal) in April, 2001. We thank Fran\c cois Lalonde for
his hospitality.


\end{document}